\documentclass[journal]{IEEEtran}

\usepackage{amsmath}
\usepackage{amssymb}
\usepackage{graphicx}
\usepackage{verbatim}
\usepackage{cite}
\usepackage{url}
\usepackage{subfigure}
\usepackage{psfrag}

\newtheorem{lemma}{Lemma}
\newtheorem{theorem}{Theorem}
\newtheorem{proposition}{Proposition}

\newcommand{\beq}{\begin{equation}}
\newcommand{\eeq}{\end{equation}}
\DeclareMathOperator{\diag}{diag}
\DeclareMathOperator{\tr}{Tr}
\DeclareMathOperator{\MSE}{MSE}
\DeclareMathOperator{\sgn}{sgn}

\newcommand{\x}{{\boldsymbol x}}
\newcommand{\y}{{\boldsymbol y}}
\newcommand{\w}{{\boldsymbol w}}
\newcommand{\z}{{\boldsymbol z}}
\renewcommand{\v}{{\boldsymbol v}}
\renewcommand{\u}{{\boldsymbol u}}
\newcommand{\Cx}{{\bf C}_\x}
\newcommand{\Cw}{{\bf C}_\w}
\newcommand{\Q}{{\bf Q}}
\newcommand{\G}{{\bf G}}
\newcommand{\D}{{\bf D}}
\renewcommand{\H}{{\bf H}}
\renewcommand{\S}{{\cal S}}
\newcommand{\T}{{\bf T}}
\newcommand{\V}{{\bf V}}
\newcommand{\I}{{\bf I}}
\newcommand{\epsmax}{\epsilon_{\max}}

\newcommand{\HCH}{\H^* \Cw^{-1} \H}
\newcommand{\hx}{{\hat{\x}}}
\newcommand{\hv}{{\hat{\v}}}
\newcommand{\hz}{{\hat{\z}}}
\newcommand{\hxls}{{\hx_{\mathrm{LS}}}}
\newcommand{\hxm}{{\hx_{\mathrm{M}}}}
\newcommand{\hxsbm}{{\hx_{\mathrm{SBM}}}}
\newcommand{\hxebm}{{\hx_{\mathrm{EBM}}}}
\newcommand{\hxbbm}{{\hx_{\mathrm{BBM}}}}
\newcommand{\hxpbm}{{\hx_{\mathrm{PBM}}}}
\newcommand{\hxt}{{\hx_{\mathrm{T}}}}
\newcommand{\E}[1]{E \kern -2pt \left\{ #1 \right\}}
\renewcommand{\Pr}[1]{\mathrm{Pr} \kern -2pt \left\{ #1 \right\}}
\newcommand{\bSigma}{\boldsymbol{\Sigma}}
\newcommand{\epsi}{\epsilon_0}
\newcommand{\lamax}{\lambda_{\max}}

\begin{document}

\title{Blind Minimax Estimation}
\author{Zvika Ben-Haim \qquad Yonina C. Eldar%
\thanks{The authors are with the Dept.\ of Electrical Engineering, Technion---Israel Institute of Technology, Haifa 32000, Israel. E-mail: zvikabh@technion.ac.il, yonina@ee.technion.ac.il; phone +972-4-829-4700; fax +972-4-829-5757. This work was supported by the Israel Science Foundation under Grant No.\ 536/04.}}
\maketitle

\begin{abstract}
We consider the linear regression problem of estimating an unknown, deterministic parameter vector based on measurements corrupted by colored Gaussian noise. We present and analyze blind minimax estimators (BMEs), which consist of a bounded parameter set minimax estimator, whose parameter set is itself estimated from measurements. Thus, one does not require any prior assumption or knowledge, and the proposed estimator can be applied to any linear regression problem. We demonstrate analytically that the BMEs strictly dominate the least-squares estimator, i.e., they achieve lower mean-squared error for any value of the parameter vector. Both Stein's estimator and its positive-part correction can be derived within the blind minimax framework. Furthermore, our approach can be readily extended to a wider class of estimation problems than Stein's estimator, which is defined only for white noise and non-transformed measurements. We show through simulations that the BMEs generally outperform previous extensions of Stein's technique.
\end{abstract}


\section{Introduction}
\label{se:intro}

The problem of estimating a parameter vector from noisy measurements has countless applications in science and engineering. Such estimation problems are typically modeled either in a Bayesian setting, in which a prior distribution on the parameter is assumed, or in a deterministic setting, in which no prior is assumed \cite{kay93}. This paper examines the deterministic estimation problem. We further assume that the measurements $\y = \H \x + \w$ are linear combinations of the parameter vector $\x$, to which Gaussian noise $\w$ is added. Here the transformation matrix $\H$ and the noise covariance are assumed to be known. We seek an estimate $\hx$ which approximates $\x$ in the sense of minimal mean-squared error (MSE).

This ubiquitous problem was first addressed by Gauss \cite{gauss1821} and Legendre \cite{legendre1806}, who proposed the classical {\it least-squares} (LS) estimator. Several lines of reasoning can be used to support the LS approach. One argument is that the LS estimator minimizes the squared error between the measurements $\y$ and the transformed estimate $\hat{\y} = \H \hx$. The LS estimator is also the maximum likelihood solution for Gaussian noise. However, neither of these criteria are directly related to the MSE, or to any other measure of the distance between $\x$ and $\hx$. Another property of the LS solution is that it is the unbiased estimator achieving minimal MSE\@. Yet by removing the requirement of unbiasedness, estimators yielding lower MSE can be constructed. While linearity and unbiasedness may be intuitively appealing properties, they have no relation to the primary goal at hand, namely, achieving low estimation error. Indeed, there are many examples in which the requirement of unbiasedness results in absurd estimators \cite{romano85}.

Because the parameter vector $\x$ is deterministic, the MSE $\E{\|\x-\hx\|^2}$ is generally a function of $\x$. In other words, one method may be better than another for some values of $\x$, and worse for other values. For instance, the trivial estimator $\hx = {\bf 0}$ achieves optimal MSE when $\x = {\bf 0}$, but its performance is otherwise poor. Nonetheless, it is possible to impose a partial order among estimation techniques \cite{lehmann98}, as follows. An estimator $\hx_1$ is said to {\it strictly dominate} a different estimator $\hx_2$ if the MSE of $\hx_1$ is lower than that of $\hx_2$, for all values of $\x$. If the MSE of $\hx_1$ is never higher than that of $\hx_2$, and is strictly lower for at least one parameter value, then $\hx_1$ is said to {\it dominate} $\hx_2$. An estimator is said to be {\it admissible} if it is not dominated by any other estimator. Surprisingly, when the parameter vector contains three or more elements, the LS method turns out to be inadmissible, i.e., some techniques {\it always} achieve lower MSE \cite{stein56}. Thus, it is of interest to characterize the class of admissible estimators, and to find techniques which dominate LS\@.

The study of admissibility is sometimes restricted to linear methods $\hx = \G \y$. A linear admissible estimator is one which is not dominated by any other linear strategy. A simple rule characterizes the class of linear admissible techniques \cite{cohen66}, and, given any linear inadmissible estimator, it is possible to construct a linear admissible alternative which dominates it \cite{eldar05}. However, the problem of admissibility is considerably more intricate when the linearity restriction is removed; generally, admissible estimators are either trivial (e.g., $\hx={\bf 0}$) or exceedingly complex \cite{maruyama98, maruyama00}. As a result, much research has focused on finding simple nonlinear techniques which dominate LS\@.

Early work on LS-dominating strategies considered the independent, identical-distribution (i.i.d.) case, for which $\H=\I$ and the noise is white. Among these, the James-Stein estimator \cite{james61,lehmann98} is the best-known example; others approaches include the works of Stein \cite{stein56} and Thompson \cite{thompson68}. Various ``extended'' James-Stein methods were later constructed for the general (non-i.i.d.) case \cite{bock75, efron73a, berger76b, manton98}. Of these, Bock's technique \cite{bock75} is quoted most often \cite{greenberg83, manton98}. However, none of these approaches has become a standard alternative to the LS estimator, and they are rarely used in practice in engineering applications \cite{manton98}. Perhaps one reason for this is that some of the estimators are poorly justified and seem counterintuitive, and as such they are sometimes regarded with skepticism (see discussion following \cite{efron73b}). Another reason is that many of these approaches (including Bock's method) result in shrinkage estimators, consisting of a gain factor multiplying the LS estimate. Shrinkage techniques can certainly be used to reduce MSE; however, in the non-i.i.d.\ case, some measurements are noisier than others, and thus a single shrinkage factor for all measurements can be considered suboptimal. Furthermore, in some applications, a gain factor has no effect on final system performance: for example, in an image reconstruction problem, multiplying the entire image by a constant does not improve quality.

In this paper, we provide a framework for generating a wide class of low-complexity, LS-dominating estimators, which are constructed from a simple, intuitive principle, called the blind minimax approach \cite{ben-haim05, ben-haim05b}. This method is used as a basis for selecting and generating techniques tailored for given problems. Many blind minimax estimators (BMEs) reduce to Stein-type methods in the i.i.d.\ case, and they continue to dominate the LS solution in the general, non-i.i.d.\ case as well. Thus, we show analytically that the proposed technique achieves lower MSE than LS, when an appropriate condition on the problem setting is satisfied. Unlike Bock's approach, BMEs may be constructed so that they are non-shrinkage, which improves their performance. Furthermore, extensive simulations show that BMEs considerably outperform Bock's method.

BMEs are based on linear minimax estimators over a bounded parameter set \cite{pinsker80, eldar04}. These are linear methods designed for a slightly different problem, in which the parameter is known to belong to a given set. The minimax approach has been thoroughly studied in this setting, and closed-form solutions are known for many types of sets. In our case, however, no prior information about the parameter set is assumed. Instead, the blind minimax approach makes use of a two-stage process (Section~\ref{se:bme}): First, a set is estimated from the measurements; next, a minimax method for this set is used to estimate the parameter itself. The result may be viewed as a simple decision rule, independent of this two-stage construction process. Indeed, our LS-dominance proofs do not rely on the method by which the techniques are generated. In particular, the dominance results do not depend on the parameter actually lying within the estimated set. Thus, the blind minimax technique provides a framework whereby many different estimators can be generated, and provides insight into the mechanism by which these techniques outperform the LS approach.

BMEs differ in the method by which the parameter set is estimated. In Section~\ref{se:sbme}, we study the case in which the estimated set is a sphere; Section~\ref{se:ebme} derives estimators based on an ellipsoidal parameter set. Section~\ref{se:js} demonstrates that several existing Stein-type methods can also be derived in the blind minimax framework. Section~\ref{se:tikhonov} compares the blind minimax approach with LS regularization techniques, while in Section~\ref{se:sim}, the BMEs are compared with other Stein-type decision rules. The paper concludes with a discussion in Section~\ref{se:discuss}.

Throughout this paper, vectors are denoted by lowercase boldface letters, and matrices by uppercase boldface letters. The $i$th component of a vector $\v$ is written as $v_i$. $\T^{1/2}$ indicates the (unique) positive semidefinite square root of a positive semidefinite matrix $\T$. The notation $\tilde{\u} \sim {\cal N}_p(\u, \Q)$ signifies that $\tilde{\u}$ is a random vector of length $p$, distributed normally with mean $\u$ and covariance $\Q$. $\|\x\|^2$ is the Euclidean norm $\x^* \x$, and $\|\x\|^2_\T$ is the $\T$-norm $\x^* \T \x$, where $\T$ is a positive definite matrix. Finally, $\diag(a_1, \ldots, a_n)$ refers to the $n \times n$ diagonal matrix whose diagonal elements are $a_1, \ldots, a_n$.

\section{Blind Minimax Estimation}
\label{se:bme}

Consider the problem of estimating an unknown deterministic parameter vector $\x \in {\mathbb C}^m$ from measurements $\y \in {\mathbb C}^n$ given by
\beq \label{eq:linear model}
\y = \H \x + \w
\eeq
where $\H \in {\mathbb C}^{n \times m}$ is a known matrix and $\w$ is a Gaussian random vector with zero mean and covariance $\Cw$. For simplicity, we assume that $\H$ is full-rank and that $\Cw$ is positive definite.

The standard solution to this regression problem is the LS approach
\beq \label{eq:def hxls}
\hxls = (\HCH)^{-1} \H^* \Cw^{-1} \y.
\eeq
The MSE of $\hxls$ does not depend on the value of $\x$, and is given by
\beq \label{eq:def epsi}
\epsi = \E{\|\hxls-\x\|^2} = \tr ( \Q^{-1} )
\eeq
where
\beq \label{eq:def Q}
\Q = \HCH.
\eeq

Despite the popularity of the LS method, other estimators are known to achieve lower MSE\@. We propose a novel strategy leading to such LS-dominating techniques, namely, the blind minimax approach. To illustrate this concept, suppose for a moment that $\x$ is known to lie within a compact parameter set $\S$. In this case, a linear minimax estimator over the set $\S$ may be constructed \cite{pinsker80, eldar04, eldar05}. This is the linear estimator $\hxm = \G \y$ minimizing the worst-case MSE among all possible values of $\x$ in $\S$,
\beq\label{eq:lin mmx}
\hxm = \arg \min_{\hx = \G\y} \max_{\x \in \S} \E{\|\hx-\x\|^2}.
\eeq
A closed form solution of \eqref{eq:lin mmx} has been previously derived for many cases of interest. Furthermore, it has been shown that any linear minimax estimator achieves lower MSE than that of the LS method, for all values of $\x$ in $\S$ \cite{ben-haim05, eldar05}. Thus, as long as {\it some} bounded set is known to contain $\x$, minimax techniques outperform the LS estimator.

BMEs utilize minimax estimators when no parameter set is known. This is done in a two-stage process:
\begin{enumerate}
\item A parameter set $\S$ is estimated from the measurements;
\item A minimax estimator designed for $\S$ is used to estimate the parameter vector $\x$.
\end{enumerate}

Various methods for estimating the parameter set $\S$ can be used, resulting in a variety of BMEs. In this paper, we consider sets of the form $\{\x: \x^* \T \x \le L^2\}$. In the next section, we examine the case $\T=\I$, in which the parameter set is spherical, resulting in a shrinkage estimator. Subsequently, in Section~\ref{se:ebme}, we discuss the more general case in which $\T = (\HCH)^b$ for some real number $b$. In both cases, closed forms are provided, and dominance over the LS method is demonstrated.

\section{The Spherical Blind Minimax Estimator}
\label{se:sbme}

In this section, we apply the blind minimax technique using a spherical parameter set $\S$ whose radius $L$ will be estimated from measurements. We assume for now that the sphere is centered on the origin, $\S = \{\x: \|\x\|^2 \le L^2 \}$. For a given value of $L$, the linear minimax estimator is \cite{eldar04}
\beq \label{eq:def sph mmx}
\hx_\text{M} = \frac{L^2}{L^2 + \epsi} \hxls,
\eeq
where $\hxls$ is the LS estimator \eqref{eq:def hxls} and $\epsi$ is the MSE \eqref{eq:def epsi} of $\hxls$. The resulting {\it spherical BME} (SBME) will have the form \eqref{eq:def sph mmx}, where $L^2$ is estimated from the measurements.

As an estimate of $L^2$, we seek a value as close as possible to $\|\x\|^2$: a smaller value would exclude the true vector $\x$ from the parameter set, while a larger value would yield an overly conservative estimator. Since $\x$ is unknown, a natural alternative is to use $\hxls$ instead. Thus, we propose to estimate $L^2$ as $\|\hxls\|^2$. Substituting into \eqref{eq:def sph mmx}, the SBME is then given by
\beq \label{eq:def sbme}
\hxsbm = \frac{\|\hxls\|^2}{\|\hxls\|^2 + \epsi} \hxls.
\eeq

In the i.i.d.\ case, the SBME reduces to the well-known Thompson estimator \cite{thompson68}. Under suitable conditions, Thompson's technique is known to strictly dominate the LS estimator, meaning that it achieves lower MSE for all values of $\x$ \cite{baranchik70}. However, the SBME is equally well-defined for the non-i.i.d.\ case. As we shall see, the SBME strictly dominates LS in the non-i.i.d.\ case, and can thus be viewed as a generalization of Thompson's results. In Section~\ref{se:js} we will demonstrate that the blind minimax approach can be used to derive generalizations of additional well-known methods, including Stein's estimator.

Up to this point, we have arbitrarily chosen the parameter set to be centered on the origin. The result was a weighted average between the LS estimate and ${\bf 0}$. Averaging with a constant value ${\bf 0}$ may be viewed as a restraint, which lessens the effect of measurement noise. As we shall see, the proposed BMEs outperform the LS estimator. This result demonstrates the fact that the LS approach results in an overestimate: reducing the norm of $\hxls$ improves its performance. However, the choice of a parameter set centered on the origin is completely arbitrary; BMEs may be constructed around any constant center point $\x_0$ \cite{greenberg83}. This will result in a weighted average between $\hxls$ and $\x_0$, which may be useful if the parameter vector is expected to lie near a particular point. Thus, the ``off-center'' SBME is given by
\beq
\hx = \left( \frac{\|\hxls\|^2}{\|\hxls\|^2+\epsi} \right) \hxls +
      \left( \frac{\epsi}{\|\hxls\|^2+\epsi} \right) \x_0.
\eeq
All dominance results continue to hold for the off-center techniques as well. In the sequel, we assume $\x_0 = {\bf 0}$ merely for the sake of notational simplicity.

The following theorem demonstrates that the SBME is guaranteed to outperform LS in terms of MSE\@.

\begin{theorem} \label{th:sbme dom lse}
Suppose $\epsi/\epsmax > 4$, where $\epsi$ is given by \eqref{eq:def epsi}, $\epsmax$ is the largest eigenvalue of $\Q^{-1}$, and $\Q=\HCH$.
Then, the SBME \eqref{eq:def sbme} strictly dominates the LS estimator.
\end{theorem}

The value $\epsi/\epsmax$ is known as the effective dimension \cite{manton98}, and may be roughly described as the number of independently-measured parameters in the system. In the i.i.d.\ case, for example, the effective dimension simply equals the length of the vector $\x$. Thus, the condition of Theorem~\ref{th:sbme dom lse} can be roughly stated as a requirement for a sufficient number of independent parameters. This requirement is a result of the fact that the LS estimator is admissible when up to two parameters are estimated \cite{stein56}. However, since many estimation problems contain dozens or hundreds of parameters and measurements, the requirement on the effective dimension holds for a variety of applications.

Note that the SBME is a special case of the estimator
\beq \label{eq:def hxb}
\hx_c = \left( 1 - \frac{\epsi}{c + \|\hxls\|^2} \right) \hxls,
\eeq
in which $c=\epsi$. Thus, rather than proving Theorem~\ref{th:sbme dom lse}, we prove the following, more general proposition, which will also be used in Section~\ref{se:js}.

\begin{proposition} \label{pr:gen sbme dom lse}
Under the conditions of Theorem~\ref{th:sbme dom lse}, the estimator $\hx_c$ given by \eqref{eq:def hxb} strictly dominates the LS estimator, for any $c \ge 0$.
\end{proposition}

The proof of Proposition~\ref{pr:gen sbme dom lse} makes use of the following lemma, which is due to Stein \cite[Theorem~1.5.15]{lehmann98}.

\begin{lemma}[Stein] \label{le:stein}
Let $\hv \sim {\cal N}_p(\v,\I)$, and let $g(\hv)$ be a differentiable function such that $\E{\left| \frac{\partial g(\hv)}{\partial \hat{v}_i} \right|} < \infty$ for all $i$.
Then,
\beq
\E{\frac{\partial g(\hv)}{\partial \hat{v}_i}} = - \E{g(\hv) (v_i - \hat{v}_i)}.
\eeq
\end{lemma}

\begin{proof}[Proof of Proposition~\ref{pr:gen sbme dom lse}]
To prove the proposition, first note that the MSE $R(\hx_c) = \E{\|\x - \hx_c\|^2}$ of $\hx_c$ is given by
\begin{multline}
R(\hx_c)
= \epsi + \E{\frac{\epsi^2 \|\hxls\|^2}{(c + \|\hxls\|^2)^2}} \\
+ 2\E{\frac{\epsi}{c+\|\hxls\|^2} \hxls^*(\x-\hxls)}. \label{eq:hbtheta_risk}
\end{multline}
Let $\V \bSigma \V^*$ be the eigenvalue decomposition of $\Q$, such that $\V$ is unitary and $\bSigma = \diag(\sigma_1, \ldots, \sigma_m)$. Define $\hv = \V^* \Q^{1/2} \hxls$ and $\v = \V^* \Q^{1/2} \x$. With these definitions, we have
\begin{align}
\hv^* \bSigma^{-1} \v  &= \hxls^* \x, \notag \\
\hv^* \bSigma^{-1} \hv &= \|\hxls\|^2, \label{eq:vDv} \\
\hv^* \bSigma^{-2} \hv &= \|\hxls\|^2_{\Q^{-1}}. \notag
\end{align}
Using these properties, the third term in \eqref{eq:hbtheta_risk} becomes
\begin{align}
&\E{\frac{\epsi}{c+\|\hxls\|^2} \hxls^*(\x-\hxls)} \notag\\
&\quad= \E{ \frac{\epsi}{c+\hv^*\bSigma^{-1}\hv} \hv^* \bSigma^{-1} (\v-\hv) } \notag\\
&\quad= \epsi \sum_{i=1}^p \sigma_i^{-1} \E{\frac{\hat{v}_i (v_i-\hat{v}_i)}{c + \hv^*\bSigma^{-1}\hv}}. \label{eq:risk3rdterm}
\end{align}

To evaluate \eqref{eq:risk3rdterm}, let
\beq
g_i(\hv) \triangleq \frac{\hat{v}_i}{c + \hv^*\bSigma^{-1}\hv},
\eeq
and note that $\hv$ is distributed normally with mean $\v$ and covariance $\I$. We can thus apply Lemma~\ref{le:stein} to obtain
\begin{align}
&\E{\frac{\epsi}{c+\|\hxls\|^2} \hxls^*(\x-\hxls)} \notag\\
&\quad= -\epsi \sum_i \sigma_i^{-1} \E{\frac{1}{c+\hv^*\bSigma^{-1}\hv} - 2 \frac{\sigma_i^{-1} \hat{v}_i^2}{(c+\hv^*\bSigma^{-1}\hv)^2}} \notag\\
&\quad= -\epsi \E{\frac{\tr(\bSigma^{-1})}{c+\hv^*\bSigma^{-1}\hv}} + 2\epsi \E{\frac{\hv^*\bSigma^{-2}\hv}{(c+\hv^*\bSigma^{-1}\hv)^2}} \notag\\
&\quad= -\epsi \E{\frac{\tr(\Q^{-1})}{c+\|\hxls\|^2}}  + 2\epsi \E{\frac{\|\hxls\|^2_{\Q^{-1}}}{(c+\|\hxls\|^2)^2}}.
\end{align}
Substituting this result back into \eqref{eq:hbtheta_risk}, we have
\begin{multline}
R(\hx_c) = \epsi +
E \Bigg\{ \frac{\epsi}{c+\|\hxls\|^2} \\
\cdot \left( \epsi \frac{\|\hxls\|^2}{c+\|\hxls\|^2}
 - 2\epsi + 4 \frac{\|\hxls\|^2_{\Q^{-1}}}{c+\|\hxls\|^2} \right)
 \Bigg\}.
\end{multline}
Since $c \ge 0$,
\begin{align}
R(\hx_c) 
&\le \epsi + \E{ \frac{\epsi}{c+\|\hxls\|^2} \left( -\epsi + 4 \epsmax \right) }.
\end{align}
If $\epsi > 4 \epsmax$, then the expectation is taken over a strictly negative range, and hence $R(\hx_c)$ is always lower than $\epsi$, so that $\hx_c$ strictly dominates $\hxls$.
\end{proof}

As we have shown, in terms of MSE, the SBME outperforms LS, providing us with a first example of the power of blind minimax estimation. The SBME is a shrinkage estimator, i.e., it consists of the LS estimator multiplied by a gain factor smaller than one. The SBME thus illustrates the fact that the LS technique tends to be an overestimate, and shrinkage can improve its performance.

\section{The Ellipsoidal Blind Minimax Estimator}
\label{se:ebme}

\subsection{Motivation}

Not all elements of the least-squares estimate $\hxls$ are equally trustworthy. Rather, $\hxls$ is a Gaussian random vector with mean $\x$ and covariance $\Q^{-1} = \left(\HCH\right)^{-1}$. Thus, some components of $\hxls$ have lower variance than others. In this sense, the scalar shrinkage factor of the SBME \eqref{eq:def sbme} and other extended Stein estimators \cite{bock75} seems inadequate.

Indeed, several researchers have proposed shrinking each measurement according to its variance. Efron and Morris \cite{efron73a} propose an empirical Bayes technique, in which high-variance components are shrunk more than low-variance ones. However, no closed form is available for this estimator, and obtaining an estimate requires iteratively solving a set of nonlinear equations. Furthermore, it is not known whether this method dominates LS\@. By contrast, Berger \cite{berger76b} provides an estimator in which more shrinkage is applied to low-variance measurements, despite the fact that low-noise components are those for which the LS approach is most accurate. Berger's technique is constructed such that the shrinkage of all components is negligible whenever there is a substantial difference between the variances of different components. As a result, dominance over the LS method is guaranteed, but the MSE gain is insubstantial unless all noise components have similar variances.

\begin{figure}
\center{\includegraphics{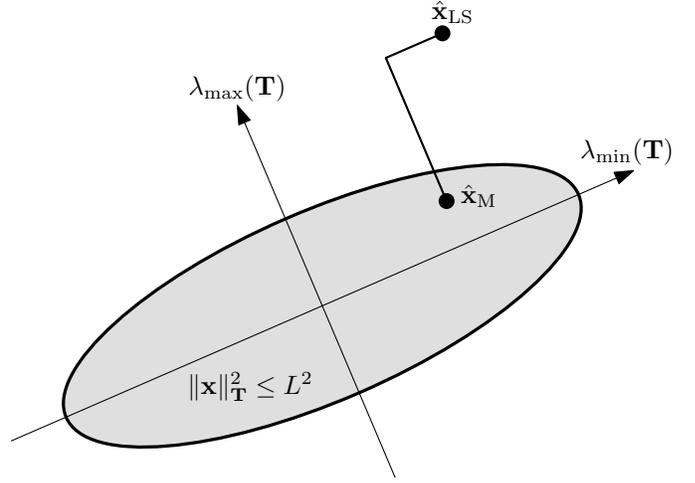}}
\caption{Illustration of the adaptive shrinkage of the minimax estimator $\hxm$ for the parameter set $\x^* \T \x \le L^2$. Low shrinkage is applied to components of $\hxls$ corresponding to small eigenvalues of $\T$, while components in directions of large eigenvalues obtain higher shrinkage.}
\label{fi:ellipse}
\end{figure}

Minimax estimators can easily be adapted for non-scalar shrinkage. Specifically, consider an ellipsoidal parameter set of the form $\S = \{ \x: \|\x\|^2_\T \le L^2 \}$, for some positive definite matrix $\T$ (see Fig.~\ref{fi:ellipse}). Let $\hxm$ represent the linear minimax estimator for this set. It can be shown that $\hxm$ is a linear function of $\hxls$, and one can therefore examine its effect on each component of $\hxls$. Consider first components of $\hxls$ in the direction of narrow axes of the ellipsoid $\S$. These components correspond to large eigenvalues of $\T$, and are denoted $\lamax(\T)$ in Fig.~\ref{fi:ellipse}. The parameter set imposes a tight constraint in these directions, and there will thus be considerable shrinkage of these elements. By contrast, components in the direction of wide axes of $\S$ (small eigenvalues of $\T$) are not constrained as tightly. Less shrinkage will be applied in this case, since the LS method is the linear minimax estimator for an unbounded set. In Fig.~\ref{fi:ellipse}, the shrinkage of wide-axis and narrow-axis components is illustrated schematically for a particular value of $\hxls$.

Typically, one would want to obtain higher shrinkage for high-variance components. Since the covariance of $\hxls$ is $\Q^{-1}$, we propose a BME based on a parameter set of the form
\beq \label{eq:ebme param set}
\S = \{ \x: \|\x\|^2_{\Q^b} \le L^2 \}
\eeq
for some constant $b<0$. The bound $L^2$ is estimated as $L^2 = \|\hxls\|^2_{\Q^b}$. We refer to the resulting technique as the {\it ellipsoidal BME} (EBME). Note that highly negative values of $b$ yield an eccentric ellipsoid, and hence result in a larger disparity between the shrinkage of different measurements. Contrariwise, a choice of $b=0$ yields scalar shrinkage, and the resulting estimator is identical to the SBME. As we will demonstrate, the EBME dominates the LS method under a condition similar to that of the SBME\@. However, the dominance condition of the EBME becomes stricter as $b$ becomes more negative. Thus, there exists a tradeoff between selective shrinkage and a broad dominance condition. In the numerical examples below we will choose a value of $b=-1$ as a compromise. 

\begin{figure*}
\centerline{%
\subfigure[]{%
\includegraphics{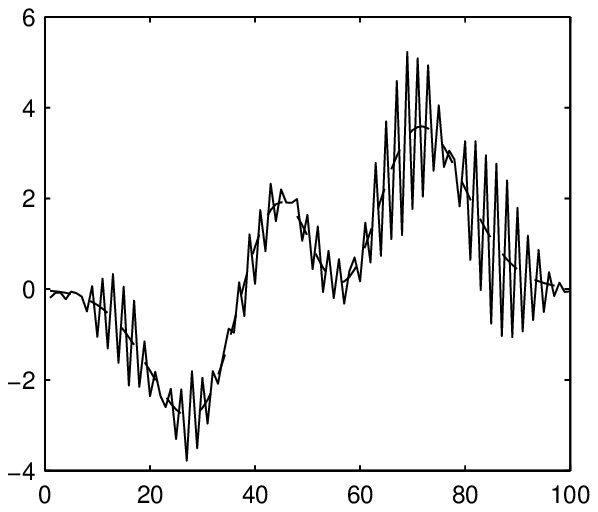}
\label{fi:example ls}}
\hfil
\subfigure[]{%
\includegraphics{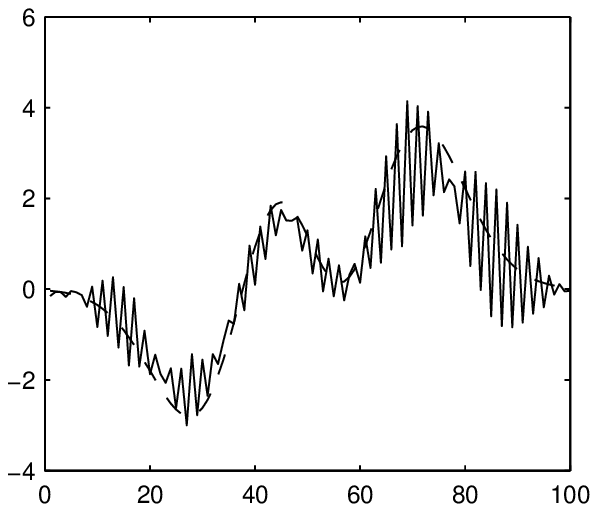}
\label{fi:example sbme}}
\hfil
\subfigure[]{%
\includegraphics{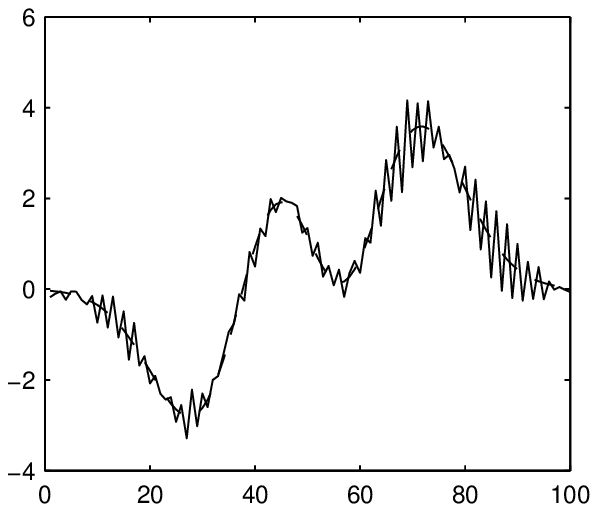}
\label{fi:example ebme}}
}
\caption{Estimation of a signal from measurements of its DCT\@. In this example, high-frequency components have a much higher noise variance than low-frequency components. Dashed line indicates original signal; solid line indicates estimate. (a) LS estimate; (b) Spherical BME, resulting in a shrinkage factor of 0.79; (c) Ellipsoidal BME, with shrinkage in the range 0.44--0.98.}
\label{fi:example}
\end{figure*}

As an additional motivation for the use of the EBME, consider the following application example (Fig.~\ref{fi:example}). Here, a 100-sample signal is to be estimated from measurements of its discrete cosine transform (DCT). Each component of the DCT is corrupted by Gaussian noise: high-variance noise is added to the 10 highest-frequency components, while the remaining components contain much lower noise levels. Thus, $\Cw$ is diagonal, and $\H$ is the DCT matrix. The condition number of $\HCH$ is 1000\@.

Since $\Cw$ is diagonal, the LS estimator is equivalent to an inverse DCT transform, and thus ignores the differences in noise level between measurements. This causes substantial estimation error, as observed in Fig.~\ref{fi:example ls}. The error is reduced by the SBME (Fig.~\ref{fi:example sbme}), which multiplies the LS estimate by an appropriately chosen scalar; in the example above, the squared error was reduced by 20\% compared with that of the LS estimate. Hence, merely multiplying the result of the LS technique by an appropriately chosen scalar can significantly reduce estimation error. However, the most significant advantage is obtained by the EBME (Fig.~\ref{fi:example ebme}), which shrinks the high-noise coefficients. Specifically, in this example, the choice $b=-1$ resulted in shrinkage of 0.44 for the high-noise coefficients, and shrinkage of only 0.98 for low-noise coefficients. The resulting squared error was 83\% lower than that of the LS estimate.

Thus, our preliminary example demonstrates that it is possible to achieve substantial improvements over the LS technique by using non-scalar shrinkage. As we will demonstrate presently, this empirical finding is only an example of the wide range of cases in which the EBME is guaranteed to improve on the LS approach.

\subsection{Dominance}

We begin our analysis by obtaining an expression for the EBMEs. A closed form solution for minimax estimators of an ellipsoidal parameter set was developed in \cite{eldar04}. By substituting the value of $L^2$ into this closed form, we obtain the following result.

\begin{proposition}[Closed-Form EBME] \label{pr:def ebme}
Let $\V \bSigma \V^*$ be the eigenvalue decomposition of $\Q = \HCH$, where $\V$ is orthonormal and $\bSigma = \diag(\sigma_1, \ldots, \sigma_m)$. Let $b \in \mathbb{R}$ be any constant, and suppose the eigenvalues $\bSigma$ are ordered such that $\sigma_1^b \ge \sigma_2^b \ge \cdots \ge \sigma_m^b > 0$. Then, the EBME for the parameter set $\S = \{ \x: \|\x\|^2_{\Q^b} \le L^2 \}$ with $L^2 = \|\hxls\|^2_{\Q^b}$ is given by
\beq \label{eq:def ebme}
\hxebm = \V \diag \left((1-\alpha \sigma_1^{b/2})_+, \ldots, (1-\alpha \sigma_m^{b/2})_+ \right) \V^* \hxls
\eeq
when $\hxls \ne {\bf 0}$, and by $\hxebm = {\bf 0}$ when $\hxls = {\bf 0}$. Here $(\cdot)_+ = \max(\cdot,0)$,
\beq
\alpha = \frac{r_1}{\|\hxls\|^2_{\Q^b} + r_2} \notag
\eeq
\beq
\begin{split}
r_1 &= \sum_{i=k+1}^m \sigma_i^{b/2-1} \label{eq:def r} \\
r_2 &= \sum_{i=k+1}^m \sigma_i^{b-1}
\end{split}
\eeq
and $k$ is chosen as the smallest index $0 \le k \le m-1$ such that
\beq \label{eq:def k}
\alpha \sigma_{k+1}^{b/2} < 1.
\eeq
\end{proposition}

\begin{proof}
In the case $\hxls={\bf 0}$, we need to find the linear minimax estimator for the set $\S = \{ {\bf 0} \}$. Clearly, the solution in this case is $\hx={\bf 0}$. For all other values of $\hxls$, we seek the linear minimax estimator for the set $\S = \{ \x : \x^* \Q^b \x \le L^2 \}$, where $L^2 = \hxls^* \Q^b \hxls > 0$. Substituting this value of $L^2$ into Proposition~1 of \cite{eldar04} yields
\begin{align} \label{eq:pr:def ebme:1}
&\hxebm = \V \diag( \underbrace{0, \ldots, 0}_{k}, \underbrace{1, \ldots, 1}_{m-k} ) \V^* (\I - \alpha \Q^{b/2}) \hxls \notag\\
&\quad  = \V \diag( \underbrace{0, \ldots, 0}_{k}, 1 - \alpha \sigma_{k+1}^{b/2}, \ldots, 1 - \alpha \sigma_m^{b/2} ) \V^* \hxls.
\end{align}
From \eqref{eq:def k}, it follows that $1 - \alpha \sigma_i^{b/2} < 0$ for all $i \le k$, and therefore \eqref{eq:pr:def ebme:1} can be written as \eqref{eq:def ebme}.
\end{proof}

We note that, as long as $\|\hxls\|^2_{\Q^b} > 0$, it is always possible to find a value $k$ which satisfies \eqref{eq:def k}. In particular, for $k=m-1$, we have
\beq
\alpha = \frac{\sigma_m^{b/2-1}}{\|\hxls\|^2_{\Q^b} + \sigma_m^{b-1}} < \frac{\sigma_m^{b/2-1}}{\sigma_m^{b-1}},
\eeq
which satisfies the requirement \eqref{eq:def k}.

While the closed form of the EBME appears somewhat more intimidating than that of the SBME, the computational complexities of the two estimators are comparable. The major difference is the calculation of the value $k$, for which $m$ divisions are required. Like the SBME, the EBME also dominates the LS estimator under suitable conditions, as shown in the following theorem.

\begin{theorem} \label{th:ebme dom lse}
Let $\hxebm$ be the EBME \eqref{eq:def ebme} and suppose that
\beq \label{eq:ebme dom cond}
\tr( \Q^{b/2-1} ) > 4 \lamax(\Q^{b/2-1})
\eeq
where $\lamax(\Q^{b/2-1})$ is the largest eigenvalue of $\Q^{b/2-1}$ and $\Q = \HCH$.
Then, $\hxebm$ strictly dominates the LS estimator.
\end{theorem}

Note that by substituting $b=0$, this result can be used to demonstrate the dominance of the SBME over LS estimation (Theorem~\ref{th:sbme dom lse}). However, the method of proof here is different, and the proof of Theorem~\ref{th:sbme dom lse} will also be used in Section~\ref{se:js}.

Also note that the dominance condition \eqref{eq:ebme dom cond} is satisfied by many reasonable estimation problems. Assuming a sufficient number of parameters, the only case in which this condition does {\it not} hold is the situation in which a small number of parameters (less than four) have much higher variance than all other parameters; in this case, the LS method is admissible or nearly so.

In order to prove Theorem~\ref{th:ebme dom lse}, we observe that the form \eqref{eq:def ebme} of the EBME is similar to Baranchik's positive-part modification \cite{baranchik64, lehmann98} of the James-Stein estimator. Baranchik proposed using a shrinkage factor of 0 whenever the James-Stein technique contains negative shrinkage, and showed that the resulting method dominates the James-Stein estimator. Although the EBME is not a shrinkage technique, it resembles Baranchik's modification, since each negative diagonal component in \eqref{eq:def ebme} is replaced with zero. The following proposition shows that the MSE can be reduced by eliminating this negative shrinkage.

\begin{proposition} \label{pr:gen pos part}
Let $\V \bSigma \V^*$ be the eigenvalue decomposition of $\Q = \HCH$, and let $b \in \mathbb{R}$ be a constant. Suppose $\hx$ is an estimator of the form $\hx = \V \D \V^* \hxls$, where $\D$ is a diagonal matrix, whose diagonal elements $d_i$ are functions of the random variable $\|\hxls\|^2_{\Q^b}$. Suppose at least one of the elements $d_i$ is negative with nonzero probability. Then, $\hx$ is dominated by the (generalized) positive-part estimator
\beq \label{eq:def gen pos part}
\hx_+ = \V \D_+ \V^* \hxls,
\eeq
where $\D_+$ is a diagonal matrix with diagonal elements $d_{i+} = \max(0,d_i)$.
\end{proposition}

\begin{proof}
Our proof follows that of Baranchik \cite{baranchik64}. We will show that $\MSE(\hx) - \MSE(\hx_+)$ is nonnegative for all $\x$, and positive for any value of $\x$ whose elements are all nonzero.
\begin{align}
\MSE&(\hx) - \MSE(\hx_+)
= \E{\| \hx-\x \|^2} - \E{\| \hx_+ - \x \|^2} \notag\\
&= \E{\|\hx\|^2 - \|\hx_+\|^2} - 2\E{\hx^*\x - \hx_+^*\x} \notag\\
&= \E{\hxls^* \V (\D^2 - \D_+^2) \V^* \hxls} \notag\\
&\quad - 2 \E{\hxls^* \V (\D - \D_+) \V^* \x}. \label{eq:pospart0}
\end{align}
Since $d_i^2 - d_{i+}^2 \ge 0$ for all $i$, the first term in \eqref{eq:pospart0} is nonnegative. Hence, to prove the proposition, it suffices to show that $\E{\hxls^* \V (\D - \D_+) \V^* \x}$ is nonpositive for all $\x$, and negative for values $\x$ with nonzero elements.

To this end, define $\z = \V^* \x$ and $\hz = \V^* \hxls$. We note that $\hz \sim {\cal N}_m(\z, \bSigma^{-1})$, so that the elements of $\hz$ are statistically independent. To calculate $\E{\hxls^* \V (\D - \D_+) \V^* \x}$, we condition on $\|\hxls\|^2_{\Q^b}$, obtaining
\begin{align}
&\E{\hxls^* \V (\D - \D_+) \V^* \x} \notag\\
&\quad= \E{ \E{\hz^* (\D-\D_+) \z | \hz^* \bSigma^b \hz } } \notag\\
&\quad= \E{ \sum_{i=1}^m (d_i - d_{i+}) \E{ \hat{z}_i z_i | \hz^* \bSigma^b \hz } } \label{eq:pospart1}
\end{align}
where we used the fact that $\|\hxls\|^2_{\Q^b} = \hz^* \bSigma^b \hz$, and that $d_i$ and $d_{i+}$ are deterministic when conditioned on $\|\hxls\|^2_{\Q^b}$. For each $i$, we further condition on $|\hat{z}_i|$, to obtain
\begin{align}
&\E{\hxls^* \V (\D - \D_+) \V^* \x} \notag\\
&\quad= \E{ \sum_{i=1}^m (d_i-d_{i+}) \E{ \hat{z}_i z_i \big| \hz^* \bSigma^b \hz, |\hat{z}_i| } } \notag\\
&\quad= \E{ \sum_{i=1}^m (d_i-d_{i+}) |\hat{z}_i z_i| \E{ \sgn(\hat{z}_i z_i) \big| \hz^* \bSigma^b \hz, |\hat{z}_i| } }.
\label{eq:pospart 2}
\end{align}
Given $|\hat{z}_i|$, we have that either $\hat{z}_i = |\hat{z}_i| \sgn(z_i)$ or that $\hat{z}_i = -|\hat{z}_i| \sgn(z_i)$. It is evident from the pdf of $\hat{z}_i$ that the latter option has lower probability, i.e.,
\begin{multline}
  \Pr{\sgn(\hat{z}_i) =   \sgn(z_i) \big| \hz^* \bSigma^b \hz, |\hat{z}_i| } \\
> \Pr{\sgn(\hat{z}_i) \ne \sgn(z_i) \big| \hz^* \bSigma^b \hz, |\hat{z}_i| }.
\end{multline}
It follows that $\E{ \sgn(\hat{z}_i z_i) | \hz^* \bSigma^b \hz, |\hat{z}_i| } \ge 0$, with strict inequality for $z_i \ne 0$. Therefore, all terms in \eqref{eq:pospart 2} are nonnegative, except for $(d_i - d_{i+})$, which is nonpositive. As a result, \eqref{eq:pospart 2} (and hence \eqref{eq:pospart0}) is nonpositive for all $\x$, so that the MSE of $\hx_+$ is never higher than that of $\hx$.

We must also show that, for some $\x$, \eqref{eq:pospart 2} is strictly negative. To this end, we choose $\x$ for which all elements are nonzero; as a result, all terms in \eqref{eq:pospart 2} are strictly positive with probability~1, except for $(d_i - d_{i+})$. The latter term is negative when $d_i<0$ and zero otherwise. Since $d_i$ is negative with nonzero probability for at least one value of $i$, we conclude that for the chosen value of $\x$, \eqref{eq:pospart 2} is strictly negative, completing the proof of Proposition~\ref{pr:gen pos part}.
\end{proof}

This generalization of the concept of a positive part estimator is now used to prove Theorem~\ref{th:ebme dom lse}.

\begin{proof}[Proof of Theorem~\ref{th:ebme dom lse}]
Clearly, the EBME \eqref{eq:def ebme} is the positive part of the estimator
\begin{align}
\hx_0 &= \V \diag \left(1-\alpha \sigma_1^{b/2}, \ldots, 1-\alpha \sigma_m^{b/2} \right) \V^* \hxls \notag\\
      &= (\I - \alpha \Q^{b/2}) \hxls.
\end{align}
Therefore, it suffices to show that $\hx_0$ dominates the LS estimator, and the theorem follows using Proposition~\ref{pr:gen pos part}.

The MSE of $\hx_0$ is given by
\begin{align} \label{eq:th:rembe 1}
&\E{ \left\| \x - \hxls + \alpha \Q^{b/2} \hxls \right\|^2 } \notag\\
&\quad= \E{ \left\| \x - \hxls + \frac{r_1 \Q^{b/2} \hxls}{\|\hxls\|^2_{\Q^b} + r_2} \right\|^2 } \notag\\
&\quad= \epsi + \E{ \frac{r_1^2 \|\hxls\|^2_{\Q^b}}{(\|\hxls\|^2_{\Q^b} + r_2)^2} } \notag\\
&\quad\quad
   + 2 \E{ \frac{r_1 (\x - \hxls)^* \Q^{b/2} \hxls}{\|\hxls\|^2_{\Q^b} + r_2} }.
\end{align}
To analyze this expression, we define
\beq \label{eq:th:gebme def v}
\begin{split}
\v  &\triangleq \V^* \Q^{1/2} \x, \\
\hv &\triangleq \V^* \Q^{1/2} \hxls.
\end{split}
\eeq
Using this notation, the third term in \eqref{eq:th:rembe 1} becomes
\begin{align} \label{eq:th:gebme 2}
A_3 &\triangleq \E{ \frac{r_1 (\x - \hxls)^* \Q^{b/2} \hxls}{\|\hxls\|^2_{\Q^b} + r_2} } \notag\\
    &= \E{ \frac{r_1 (\x - \hxls)^* \Q^{1/2} \V \V^* \Q^{b/2-1} \V \V^* \Q^{1/2} \hxls}
                {\hxls^* \Q^{1/2} \V \V^* \Q^{b-1} \V \V^* \Q^{1/2} \hxls + r_2} } \notag\\
    &= \E{ \frac{r_1 (\v - \hv)^* \bSigma^{b/2-1} \hv}{\hv^* \bSigma^{b-1} \hv + r_2} } \notag\\
    &= \sum_{i=1}^m \sigma_i^{b/2-1} \E{ \frac{ r_1 (v_i - \hat{v}_i)\hat{v}_i }{ \hv^* \bSigma^{b-1} \hv + r_2 } }.
\end{align}

Next, define
\beq \label{eq:th:gebme def g_i}
g_i(\hv) \triangleq \frac{r_1 \hat{v}_i}{\hv^* \bSigma^{b-1} \hv + r_2}.
\eeq
Note that $r_1$ and $r_2$ are implicitly dependent on $k$, which in turn depends on $\hv$. Thus, $g_i(\hv)$ is discontinuous for some values of $\hv$, namely, those values for which $\alpha = \sigma_i^{b/2}$. However, these values of $\hv$ occur with probability zero; for all other values, $k$ (and hence $r_1$ and $r_2$) are constant for sufficiently small changes in $\hv$. Thus,
\beq \label{eq:th:gebme deriv g_i}
\frac{\partial g_i}{\partial \hat{v}_i} = \frac{r_1}{\hv^* \bSigma^{b-1} \hv + r_2}
    - \frac{ 2 r_1 \sigma_i^{b-1} \hat{v}_i^2 }{ ( \hv^* \bSigma^{b-1} \hv + r_2 )^2 }
\quad \text{w.p.~1}
\eeq
and $\E{ | \partial g_i / \partial \hat{v}_j | } < \infty$ for all $i,j$. Furthermore, observe that $\hv \sim {\cal N}_m(\v,\I)$. We can therefore apply Lemma~\ref{le:stein} to $g_i$. This yields
\begin{multline}
\E{ \frac{ r_1 \hat{v}_i (v_i - \hat{v}_i) }{ \hv^* \bSigma^{b-1} \hv + r_2 } }\\
= \E{ - \frac{ r_1 }{ \hv^* \bSigma^{b-1} \hv + r_2 } + 2 \frac{ r_1 \sigma_i^{b-1} \hat{v}_i^2 }{( \hv^* \bSigma^{b-1} \hv + r_2 )^2} }.
\end{multline}
Substituting into \eqref{eq:th:gebme 2}, we obtain
\begin{align}
A_3 &= \sum_{i=1}^m \sigma_i^{b/2-1} \E{ \frac{ 2r_1 \sigma_i^{b-1} \hat{v}_i^2 }{( \hv^* \bSigma^{b-1} \hv + r_2 )^2}
                                      - \frac{ r_1 }{ \hv^* \bSigma^{b-1} \hv + r_2 } } \notag\\
    &= \E{ \frac{ 2r_1 \sum_{i=1}^m \sigma_i^{3b/2-2} \hat{v}_i^2 }{( \hv^* \bSigma^{b-1} \hv + r_2 )^2}
         - \frac{ r_1 \sum_{i=1}^m \sigma_i^{b/2-1} }{ \hv^* \bSigma^{b-1} \hv + r_2 } } \notag\\
    &= \E{ \frac{ 2r_1 \hv^* \bSigma^{3b/2-2} \hv }{( \hv^* \bSigma^{b-1} \hv + r_2 )^2}
         - \frac{ r_1 \tr(\bSigma^{b/2-1}) }{ \hv^* \bSigma^{b-1} \hv + r_2 } } .
\end{align}
Using the definition \eqref{eq:th:gebme def v} of $\hv$, $A_3$ may be written as
\beq \label{eq:th:gebme 3}
\E{ \frac{ r_1 }{ \|\hxls\|^2_{\Q^b} + r_2 }
\left[ \frac{2 \|\hxls\|^2_{\Q^{3b/2-1}}}{\|\hxls\|^2_{\Q^b} + r_2} - \tr(\Q^{b/2-1}) \right] }.
\eeq
Note that
\begin{align}
\frac{\|\hxls\|^2_{\Q^{3b/2-1}}}{\|\hxls\|^2_{\Q^b} + r_2}
    &< \frac{\|\hxls\|^2_{\Q^{3b/2-1}}}{\|\hxls\|^2_{\Q^b}} \notag\\
    &= \frac{(\Q^{b/2} \hxls)^* \Q^{b/2-1} (\Q^{b/2} \hxls)}{(\Q^{b/2} \hxls)^* (\Q^{b/2} \hxls)} \notag\\
    &\le \lamax(\Q^{b/2-1}).
\end{align}
Thus
\beq
A_3 < \E{ \frac{ r_1 }{ \|\hxls\|^2_{\Q^b} + r_2 } \left[ 2 \lamax(\Q^{b/2-1}) - \tr(\Q^{b/2-1}) \right] }.
\eeq

Substituting back into \eqref{eq:th:rembe 1}, we have
\begin{multline}
\MSE < \epsi + E\Bigg\{ \frac{r_1}{\|\hxls\|^2_{\Q^b} + r_2}\\
         \cdot \left[ r_1 + 4 \lamax(\Q^{b/2-1}) - 2 \tr(\Q^{b/2-1}) \right] \Bigg\}
\end{multline}
and using the fact that $r_1 \le \tr(\bSigma^{b/2-1}) = \tr(\Q^{b/2-1})$, we conclude that the MSE is bounded by
\beq
\epsi + \E{ \frac{r_1}{\|\hxls\|^2_{\Q^b} + r_2} \left[ 4 \lamax(\Q^{b/2-1}) - \tr(\Q^{b/2-1}) \right] }.
\eeq
Thus, if $\tr(\Q^{b/2-1}) > 4 \lamax(\Q^{b/2-1})$, then $\MSE < \epsi$, proving that the EBME dominates the LS estimator.
\end{proof}

Thus far, we have presented two examples of BMEs which dominate the LS method under suitable conditions. Both approaches are extensions of Thompson's technique to the non-i.i.d.\ case. In the next section, we demonstrate that other BMEs extend different LS-dominating techniques, namely Stein's estimator and Baranchik's positive-part improvement.

\section{Relation to Stein-type Estimation}
\label{se:js}

In Section~\ref{se:sbme}, the SBME \eqref{eq:def sbme} was constructed by using $L^2 = \|\hxls\|^2$ as an estimate of $\|\x\|^2$. However, the fact that shrinkage techniques such as the SBME dominate LS indicates that $\hxls$ is in fact an overestimate of $\x$. It is arguably more accurate to use a smaller value than $\|\hxls\|^2$ to estimate $\|\x\|^2$. In particular, it is readily shown that
\beq
\E{\|\hxls\|^2} = \|\x\|^2 + \epsi.
\eeq
Hence, one may opt to use
\beq \label{eq:L2 balanced}
L^2 = \|\hxls\|^2-\epsi
\eeq
as an estimate of $\|\x\|^2$. It is important to note that such a value of $L^2$ cannot be used with the linear minimax method, since $L^2$ is negative with nonzero probability; a parameter set with negative radius is undefined. However, substituting \eqref{eq:L2 balanced} into a minimax technique, as per the blind minimax approach, can still lead to well-defined estimators. In particular, substituting \eqref{eq:L2 balanced} into the spherical minimax method \eqref{eq:def sph mmx} yields the ``balanced'' BME
\beq \label{eq:def balanced bme}
\hxbbm = \left( 1 - \frac{\epsi}{\|\hxls\|^2} \right) \hxls.
\eeq

A striking property of the balanced BME is that it reduces to Stein's estimator \cite{stein56} in the i.i.d.\ case. Both techniques are well-defined unless $\hxls = {\bf 0}$, an event which has zero probability. Furthermore, the balanced BME extends Stein's method, in that it continues to dominate LS for the non-i.i.d.\ case, under suitable conditions. This is shown by the following theorem.

\begin{theorem} \label{th:bbme dom lse}
Suppose $\epsi/\epsmax > 4$, where $\epsi$ is given by \eqref{eq:def epsi}, $\epsmax$ is the largest eigenvalue of $\Q^{-1}$, and $\Q$ is given by \eqref{eq:def Q}.
Then, the balanced BME \eqref{eq:def balanced bme} strictly dominates the LS estimator.
\end{theorem}

\begin{proof}
The theorem follows by substituting $c=0$ in Proposition~\ref{pr:gen sbme dom lse}.
\end{proof}

A well-known drawback of Stein's approach is that it sometimes causes negative shrinkage, i.e., the shrinkage factor in \eqref{eq:def balanced bme} is negative with nonzero probability. This is known to increase the MSE \cite{baranchik64}. From the blind minimax perspective, this negative shrinkage is a result of the fact that $L^2$ can become negative. Thus, it is natural to replace \eqref{eq:L2 balanced} with
\beq \label{eq:L2 positive}
L^2 = \left( \|\hxls\|^2-\epsi \right)_+
\eeq
where $(a)_+ = \max(a,0)$. Substituting this value of $L^2$ into the spherical minimax estimator yields the ``positive-part BME,'' given by
\beq \label{eq:def pp bme}
\hxpbm = \left( \frac{ ( \|\hxls\|^2 - \epsi )_+ }{ ( \|\hxls\|^2 - \epsi )_+ + \epsi } \right) \hxls.
\eeq
Note that when $\|\hxls\|^2 - \epsi < 0$, the estimator $\hxpbm$ equals ${\bf 0}$; in all other cases, $\hxpbm = \hxbbm$. Thus, \eqref{eq:def pp bme} may be written as
\beq \label{eq:pp bme 2}
\hxpbm = \left( 1 - \frac{\epsi}{\|\hxls\|^2} \right)_+ \hxls.
\eeq
In other words, $\hxpbm$ is the positive part of the balanced BME\@. Specifically, in the i.i.d.\ case, $\hxpbm$ is the positive-part correction of Stein's estimator. In the i.i.d.\ case, Baranchik \cite{baranchik64} demonstrated that $\hxpbm$ dominates $\hxbbm$. An interesting question for further research is whether the dominance property holds in the non-i.i.d.\ case as well.

\begin{figure}
\center{\includegraphics{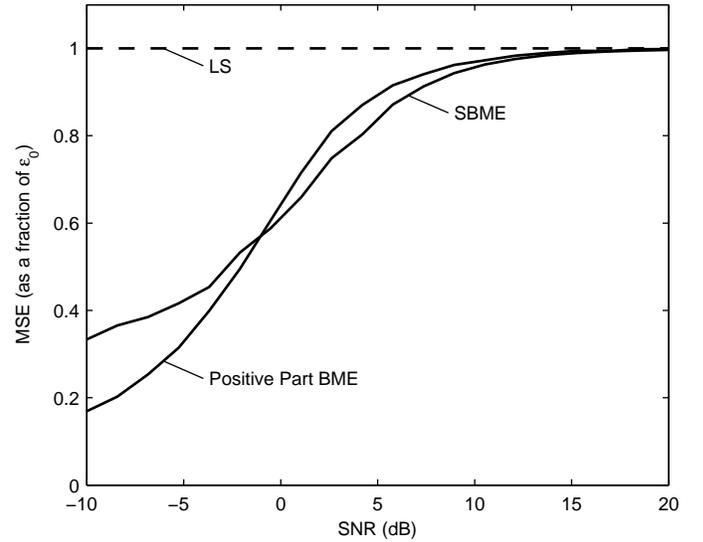}}
\caption{Comparison between the positive part approach and the SBME\@. The positive part method results in stronger shrinkage, which improves performance for low SNR at the expense of high SNR\@.}
\label{fi:positive part}
\end{figure}

The ``balanced'' method presented in this section for estimating the parameter set radius results in a value \eqref{eq:L2 balanced} of $L^2$ which is smaller than that of the SBME\@. As a result, the balanced approach causes more shrinkage towards the origin. This tends to improve performance for low signal-to-noise ratio (SNR) at the expense of performance degradation for high SNR\@. In particular, $\hxpbm$ has a positive probability of yielding an estimate of ${\bf 0}$. This may indeed reduce the MSE when the parameter is exceedingly small with respect to the noise variance, but will sacrifice high-SNR performance.

In Fig.~\ref{fi:positive part}, the positive part estimator $\hxpbm$ is compared with the SBME of Section~\ref{se:sbme}. The problem setting of this simulation is identical to that of Fig.~\ref{fi:mseplot max}, which will be described in detail in Section~\ref{se:sim}. In general, the positive-part BME tends to perform as well or worse than the SBME at SNR values above 0~dB, and better for lower SNR values. Thus, in most applications, use of the SBME is probably preferable. However, the fact that Stein's estimator can be derived and extended using blind minimax considerations illustrates the versatility of this approach.

\section{Comparison with LS Regularization}
\label{se:tikhonov}

Independently of the development of Stein-type estimators, many researchers became aware of deficiencies of the LS approach for solving ill-conditioned problems. A variety of alternatives were proposed as a result. These substitutes were generally not required to dominate the LS estimator; rather, they were intended to improve estimation quality in specific scenarios. Of these approaches, the most common is Tikhonov regularization \cite{tikhonov77}, also referred to as ridge regression \cite{hoerl70}.

Tikhonov regularization is intended for ill-posed problems, i.e., problems in which $\HCH$ is nearly singular. The matrix $\Q = \HCH$ is guaranteed to be positive-definite (and hence invertible), since we assume that $\H$ is full-rank and $\Cw$ is positive-definite. However, $\Q$ may contain eigenvalues which are very close to zero. In these cases, the LS estimator (which depends on the term $\Q^{-1}$) causes severe amplification of measurement noise. In effect, an ill-posed setting is one in which the SNR of at least one parameter is extremely low; as we have seen, the LS approach results in overestimation in such conditions. Regularization techniques attempt to mitigate this problem by improving the conditioning of the matrix $\Q$.

Tikhonov regularization may be justified in a Bayesian setting, as follows. Suppose that the parameter vector $\x$ is known to be distributed normally, independently of the noise $\w$, with zero mean and a covariance matrix $\Cx$. The minimum MSE estimator of $\x$ given $\y$ is then the Wiener filter \cite{wiener49,kay93}
\beq \label{eq:def gen tikhonov}
\hx = \bigl( \HCH + \Cx^{-1} \bigr)^{-1} \H^* \Cw^{-1} \y.
\eeq
In practice, $\x$ is a deterministic parameter, and thus does not have a covariance matrix. However, by replacing $\Cx^{-1}$ with an appropriately chosen regularization matrix, the (generalized) Tikhonov estimator is obtained.

There are several methods for empirically selecting a regularization matrix $\Cx^{-1}$. If nothing is known about the parameter $\x$, one possibility is to choose $\Cx = \sigma_x^2 \I$, where $\sigma_x^2$ is to be estimated from $\y$. Optimally, one would like to use the average value of $x_i^2$ as an approximation of the variance $\sigma_x^2$. However, since $\x$ is unknown, this is not possible. Instead, $\sigma_x^2$ can be estimated as $\sum \hat{x}_{\text{LS}, i}^2 / m$, which is an approximation of the desired quantity $\sum x_i^2 / m$. This results in the estimator
\beq \label{eq:def hxt1}
\hxt^{(1)} = \left( \HCH + \frac{m}{\|\hxls\|^2} \I \right)^{-1} \H^* \Cw^{-1} \y.
\eeq
This derivation is based on an empirical Bayes approach \cite{robbins64}, in which the elements of $\x$ are assumed to be i.i.d\@. An alternative is to assume instead that the variance of $\x$ is proportional to the variance of the noise $\w$, which implies $\Cx = \alpha \Q^{-1}$. In analogy to the previous derivation, one may then estimate $\alpha$ as $m / \|\hxls\|^2_\Q$. Substituting into \eqref{eq:def gen tikhonov} results in the shrinkage estimator
\beq \label{eq:def hxt2}
\hxt^{(2)} = \frac{\|\hxls\|^2_\Q}{m + \|\hxls\|^2_\Q} \hxls.
\eeq

\begin{figure}
\psfrag{Tikh1}{{\footnotesize $\hxt^{(1)}$}}\psfrag{Tikh2}{{\footnotesize $\hxt^{(2)}$}}
\center{\includegraphics{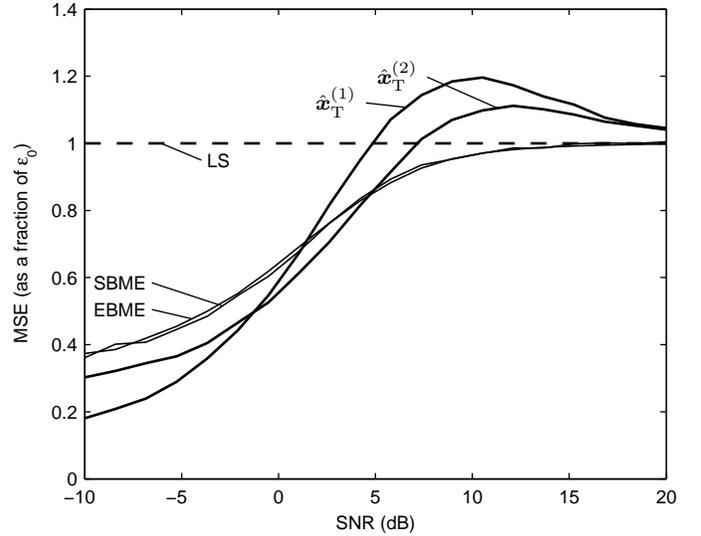}}
\caption{Tikhonov regularization does not dominate the LS estimator. The Tikhonov estimators $\hxt$ are seen to perform worse than the LS estimator at high SNR, whereas the BMEs dominate the LS method.}
\label{fi:tikhonov}
\end{figure}

Unfortunately, the Tikhonov estimators $\hxt^{(1)}$ and $\hxt^{(2)}$ do not dominate LS; like the original Tikhonov regularization, they perform poorly at high SNR values. To illustrate this, we performed a simulation in which the MSE of the LS method was compared to that of $\hxt^{(1)}$ and $\hxt^{(2)}$. In this example, 15 parameters were estimated using 15 independent measurements, with $\H=\I$. The noise variance of five of the measurements was 100 times larger than the noise variance of the remaining measurements. The parameter vector was chosen in the direction of a high-variance measurement, and its magnitude was varied to obtain different SNR values. Here and in the remainder of the paper, we define the SNR as
\beq \label{eq:def snr}
\text{SNR} = \frac{\|\x\|^2}{\E{\|\w\|^2}} = \frac{\|\x\|^2}{\tr(\Cw)}.
\eeq
For comparison, the MSE of the LS and blind minimax techniques were also calculated.

The results are displayed in Fig.~\ref{fi:tikhonov}. It is evident from this figure that the Tikhonov regularization is inadequate at high SNR, as it performs worse than the LS estimator. Both Tikhonov approaches converge to the LS approach at infinite SNR, but consistently obtain higher MSE than the LS method for SNR values above 5~dB\@. This makes them unattractive candidates for replacing the LS technique.

\section{Numerical Results}
\label{se:sim}

Estimator performance depends on a variety of operating conditions, including the effective dimension, the SNR, the eigenvalues of $\Q=\HCH$, and the value of the unknown parameter vector $\x$. Several computer simulations were implemented to test the effect of these conditions on performance of the SBME and EBME\@. In these tests, a value of $b=-1$ was used for the parameter set \eqref{eq:ebme param set} of the EBME\@. The simulations were also used to compare the BMEs with Bock's estimator \cite{bock75}, which is the most commonly-used extended Stein estimator \cite{greenberg83,manton98}. Like Stein's results, Bock's approach consists of a shrinkage estimator, given by
\beq \label{eq:def bock}
\hx_{\text{Bock}} = \left( 1 - \frac{\epsi/\epsmax-2}{\|\hxls\|^2_\Q} \right) \hxls.
\eeq

\begin{figure*}
\centerline{%
\subfigure[]{%
\includegraphics{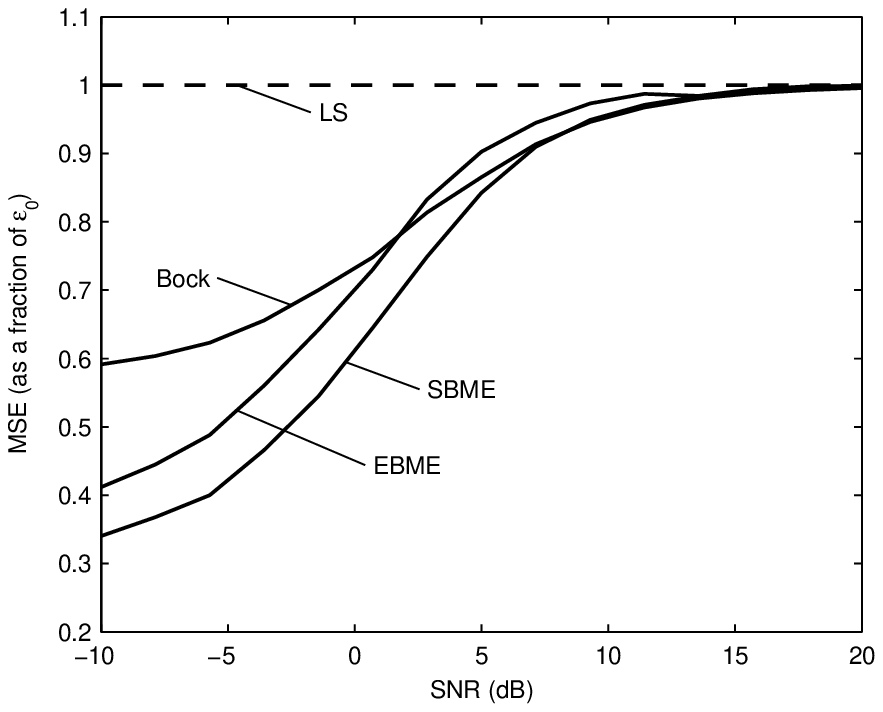}
\label{fi:mseplot max}}
\hfil
\subfigure[]{%
\includegraphics{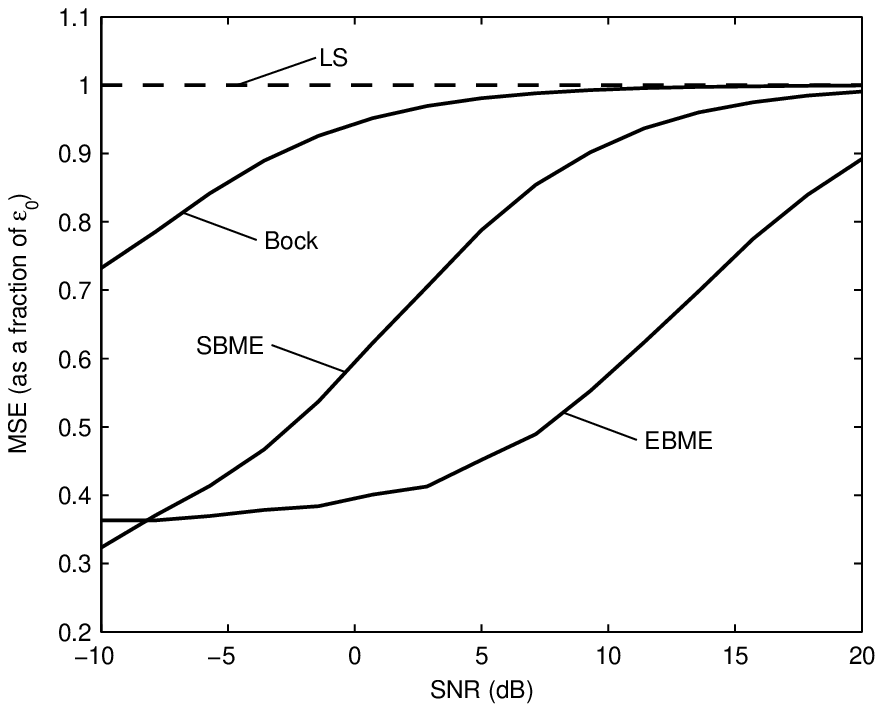}
\label{fi:mseplot min}}
}
\caption{MSE vs.\ SNR for a typical operating condition: effective dimension 5.8, $m=n=15$. (a)~Parameter vector $\x$ in direction of maximum noise; (b)~Parameter vector $\x$ in direction of minimum noise.}
\label{fi:mseplot}
\end{figure*}

The theorems of Sections \ref{se:sbme} and \ref{se:ebme} ensure that the BMEs achieve lower MSE than the LS estimator, but do not guarantee that this improvement is substantial. To measure this performance gain, we first chose a typical scenario, in which the number of parameters $m$ and the number of measurements $n$ were both 15. The system matrix $\H$ was chosen as $\I$, and the noise covariance $\Cw$ was
\beq
\Cw = \sigma^2 \diag(1,1,1,1,.5,.2,.2,.2,.2,.1,.1,.1,.1,.05,.05)
\eeq
resulting in an effective dimension of 5.8. Here $\sigma^2$ was selected to achieve the desired SNR \eqref{eq:def snr}. To illustrate the dependence on the value of the parameter vector $\x$, two different settings were tested. In Fig.~\ref{fi:mseplot max}, $\x$ is chosen in the direction of the maximum eigenvector of $\Q^{-1}$, while in Fig.~\ref{fi:mseplot min}, $\x$ is chosen in the direction of the minimum eigenvector. This corresponds to parameters in the direction of maximal and minimal noise, respectively. Estimates of the MSE were calculated for a range of SNR values by generating 10,000 random realizations of noise per SNR value.

It is evident from Fig.~\ref{fi:mseplot} that substantial improvement in MSE can be achieved by using BMEs in place of the LS approach: in some cases the MSE of the LS estimator is nearly three times larger than that of the BMEs. The performance gain is particularly noticeable at low and moderate SNR\@. At infinite SNR, the LS technique is known to be optimal \cite{kay93}, and all other methods converge to the value of the LS estimate; as a result, performance gain is smaller at high SNR, although substantial improvement can be obtained even at 10--15~dB\@.

\begin{figure*}
\centerline{%
\subfigure[]{%
\includegraphics{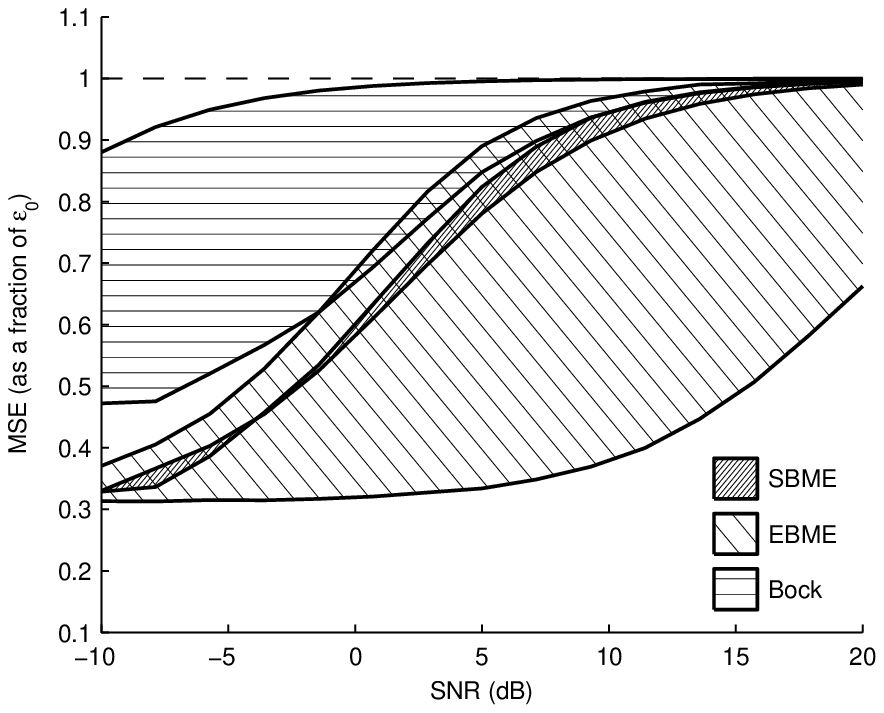}
\label{fi:1plot-lin}}
\hfil
\subfigure[]{%
\includegraphics{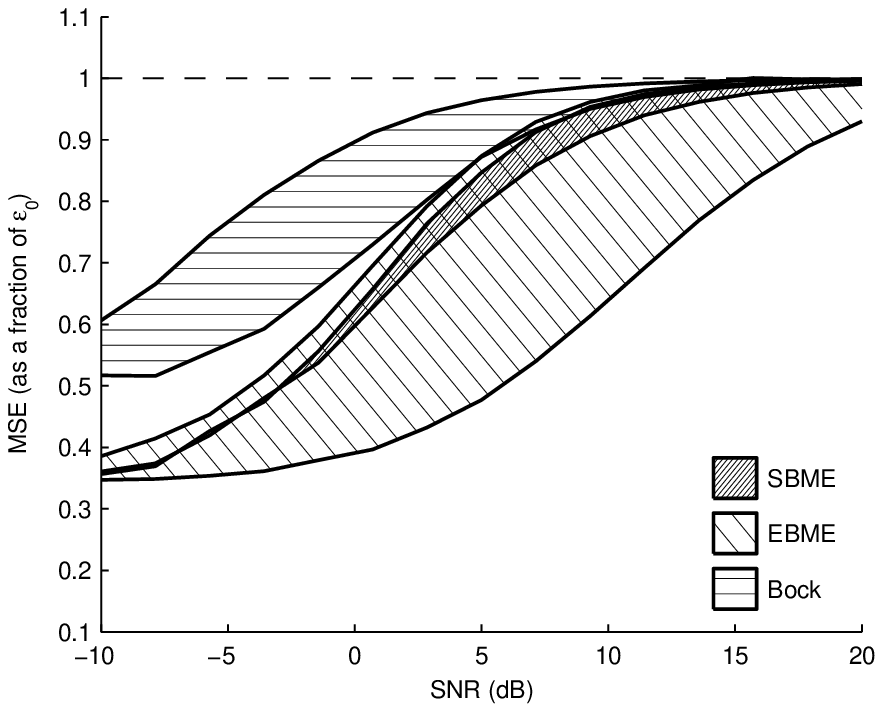}
\label{fi:1plot-stag}}
}
\caption{Range of possible MSE values obtained for different values of $\x$, as a function of SNR\@. $\H=\I$ for both figures. (a) $m=n=15$, with eigenvalues of $\Cw$ distributed uniformly between 1 and 0.01, resulting in an effective dimension of 7.6; (b) $m=n=10$, with $\Cw$ containing five eigenvalues of 1 and five eigenvalues of 0.1, resulting in an effective dimension of 5.5\@.}
\label{fi:1plot}
\end{figure*}

To further compare the BMEs with Bock's estimator, another simulation was performed, in which a large set of parameter values $\x$ were generated for different SNRs. For each estimator, and for each SNR, the lowest and highest MSE were determined, resulting in a measure of the performance range for each estimator. This performance range is displayed in Fig.~\ref{fi:1plot} for two different choices of $\Cw$, which are indicated in the figure caption. One may observe that both BMEs outperform Bock's estimator under nearly all circumstances. It is also interesting to note that while the MSE of the EBME is highly dependent on the value of the parameter value $\x$, the performance of the SBME is fairly constant. This is a result of the symmetric form of the SBME\@. On the other hand, the EBME achieves considerably lower MSE for most values of the parameter vector.

It is insightful to compare the performance of the SBME and EBME in Figs. \ref{fi:mseplot} and \ref{fi:1plot}. While the worst-case performance of the two blind minimax techniques is similar, the EBME performs considerably better for some values of $\x$. This is a result of the fact that the EBME selectively shrinks the noisy measurements, whereas the SBME uses an identical shrinkage factor for all elements. If one measurement contains very little noise, the SBME is forced to reduce the shrinkage of all other measurements. The EBME, by contrast, can effectively reduce the effect of noisy measurements without shrinking the clean elements. As a result, the EBME is superior by far if $\x$ is orthogonal to the noisiest measurements, whence the selective shrinkage is most effective; its performance gain is less substantial when $\x$ is in the direction of high shrinkage, since in these cases, shrinkage is applied to the parameter as well as the noise.

\begin{figure}
\center{\includegraphics{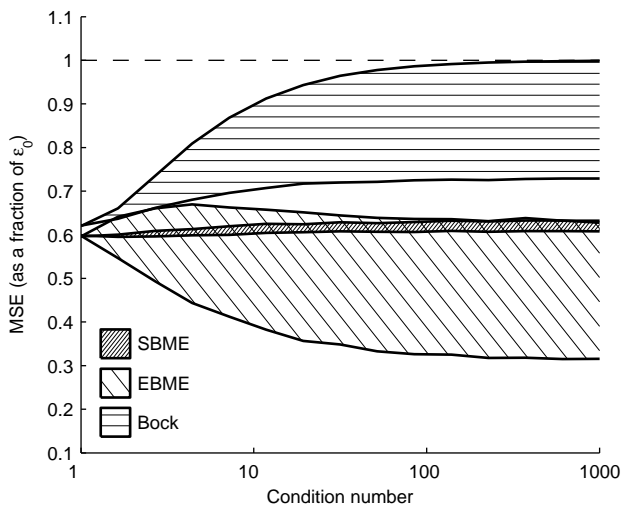}}
\caption{Range of possible MSE values obtained for different values of $\x$, as a function of the condition number of $\Q$. SNR 0~dB, $m=n=10$.}
\label{fi:cn}
\end{figure}

Another important advantage of the blind minimax approach over Bock's estimator is that the latter converges to the LS technique when the matrix $\Q$ is ill-conditioned, i.e., when some eigenvalues are much larger than others. This is because the shrinkage in Bock's method \eqref{eq:def bock} is a function of $1/\|\hxls\|^2_\Q$. As a result, when $\hxls$ contains a significant component in the direction of a large eigenvalue of $\Q$, shrinkage becomes negligible. Yet, in this case, shrinkage is still desirable for the remaining eigenvalues. This effect is demonstrated in Fig.~\ref{fi:cn}, which plots the performance of the various approaches for matrices $\Q$ having condition numbers between 1 and 1000\@. Here, 10 parameters and 10 measurements are used, $\H=\I$, and the noise covariance is chosen such that the first five eigenvalues equal 1 and the remaining five eigenvalues equal a value $v$, which is chosen to obtain the desired condition number.
For each condition number, a large set of values $\x$ are chosen such that the SNR is 0~dB; as in Fig.~\ref{fi:1plot}, the range of MSE values obtained for each estimate is plotted.
It is evident that Bock's estimator approaches the LS method for ill-conditioned matrices, despite the fact that shrinkage can still improve performance, as indicated by the performance of the SBME\@. The performance of the EBME improves relative to the LS estimator for ill-conditioned matrices, since the high-noise components are further reduced in this case.

\section{Discussion}
\label{se:discuss}

The blind minimax approach is a general technique for using minimax estimators in situations for which no parameter set is known. We considered an application of this concept to the Gaussian linear regression model. Two novel estimators were proposed: a technique based on a spherical parameter set, and one based on an ellipsoidal parameter set. In Sections \ref{se:sbme} and \ref{se:ebme}, these approaches were shown to dominate the LS method. Under fairly weak conditions, in any application which makes use of the LS estimator, the MSE performance can be improved by using a BME instead. Furthermore, in Section~\ref{se:js}, we demonstrated that Stein's approach, as well as its positive part modification, can be derived and generalized using the blind minimax framework.

It can readily be shown that the dominance condition of the SBME (Theorem~\ref{th:sbme dom lse}) is weaker than the dominance condition of the EBME (Theorem~\ref{th:ebme dom lse}), i.e., the conditions for SBME dominance hold whenever the conditions for EBME dominance hold. The dominance condition of Bock's estimator \cite{bock75} is still weaker\footnote{A simple change to the SBME (adding $-2$ to the numerator) changes its dominance condition to that of Bock's estimator, without significantly affecting its performance. However, we have been unable to derive this modification using the blind minimax approach, and thus prefer the simpler form of the SBME used in the paper.}. This would seem to indicate that Bock's estimator is superior to the proposed estimators. Yet the results of Section~\ref{se:sim} demonstrate that the opposite is true: the BMEs usually outperform Bock's estimator. This is true in particular for ill-conditioned problems, for which the LS estimator is notoriously inaccurate; for such problems, Bock's approach dominates the LS method by a negligible margin, whereas the BMEs achieve a significant performance gain. Thus, while dominance theorems are useful in providing sufficient conditions for improving on the LS estimator, they are ill-suited for comparing LS-dominating estimators. This conclusion is noteworthy since estimators are sometimes chosen by maximizing the range of conditions for which dominance is guaranteed. It seems that other analytical tools are required for comparing LS-dominating estimators. For example, it may be possible to prove that BMEs dominate Bock's estimator, for some problem settings.

The choice between the different BMEs is application-dependent. As demonstrated in Section~\ref{se:sim}, the SBME reliably achieves constant performance for a variety of values of $\x$, although the typical performance of the EBME is superior. The EBME is particularly well-adapted to ill-posed problems, in which some measurements are much more noisy than others. In such cases, the use of a single shrinkage factor for all measurements is clearly suboptimal. As a result, scalar shrinkage methods such as the SBME and Bock's technique often result in little improvement over the LS estimator, while the EBME is capable of selectively shrinking the noisy measurements, thus improving performance.

The use of a componentwise shrinkage technique such as the EBME may be useful in additional contexts as well. In some applications, MSE minimization is only a nominal goal which approximates some other error criterion. In these cases, a shrinkage estimator has no impact on the actual objective. For example, if the vector $\x$ is an image which is to be reconstructed, its subjective quality is not affected by multiplying the entire estimate by a scalar. Likewise, in a binary receiver, the sign of $\x$ must be determined, but the sign does not change when the estimate is shrunk. In such applications, the SBME (and Bock's estimator) have no effect on the final result, whereas the EBME can be used to improve performance.

\section{Conclusion}
\label{se:conclu}

In this paper, we presented the blind minimax strategy, whereby one uses linear minimax estimators whose parameter set is itself estimated from measurements. This simple concept was examined in the setting of a linear system of measurements with colored Gaussian noise, where we have shown that the BMEs dominate the LS method. Hence, in any such problem, the proposed estimators can be used in place of the LS approach, with a guaranteed performance gain. Apart from being useful in and of themselves, the proposed techniques support the underlying concept of blind minimax estimation. This concept can be applied to many other problems, such as estimation with uncertain system matrices, estimation with non-Gaussian noise, and sequential estimation. Use of the blind minimax approach in such problems remains a topic for further study.

Stein's discovery of LS-dominating estimators, half a century ago, shocked the statistics community, and his results are still rarely used in practice. It is our hope that the blind minimax concept will provide additional support for such estimators, both by supplying an intuitive understanding of Stein's phenomenon, and by providing a wide class of powerful new estimators.

\section{Acknowledgement}
\label{se:ack}

The authors are grateful to the anonymous reviewers for their careful reading of the paper, which helped clarify and improve several results.


\end{document}